\documentclass[a4paper, 12pt,oneside]{article}

\usepackage{amsthm,amssymb,amsmath}

\usepackage{hyperref}
\hypersetup{hidelinks}

\newtheorem{theorem}{Theorem}[section]

\newtheorem{lemma}[theorem]{Lemma}

\numberwithin{equation}{section}

\begin{document}
\begin{center}
{\Large \bf $L^{\Vec{p}}-L^{\Vec{q}}$ Boundedness of Multiparameter Forelli-Rudin Type Operators on the Siegel Upper Half-space }
\end{center}

\begin{center}
{\bf Hongheng Yin, Guan-Tie Deng and Zhi-Qiang Gao }
\end{center}

{\bf Abstract:}
In this article,we present exactly when two classes of multiparameter Forelli-Rudin type integral operators are bounded from one weighted mixed-norm Lebesgue space $L^{\Vec{p}}$ to another space $L^{\Vec{q}}$ over the Siegel upper half-space.


\section{Introduction}
\label{sec:int}

The main motivation for this paper is the work of Huang,Wang,Zhang\cite{BondYur18},the authors investigated the boundedness of two classes of multiparameter intergral operators induced by Bergman type kernels on weighted mixed-norm Lebesgue space on the unit ball in $\mathbb C^n$. They also provided a detailed introduction on the research history of mixed-norm Lebesgue spaces and related research on the boundedness of integral operators on the unit ball.

In this article, we will replace the unit ball with Siegel upper half-space to continue studying the boundedness of operators. Before starting our work, we will provide some basic definitions and notation, as well as previous works.

Let n be a positive integer,$\mathbb C^n=\mathbb C\times\dots \times \mathbb C$ denote the n-dimensional complex Euclidean space.The Siegel upper half-space in $\mathbb C^n$ is the set
$$
\mathcal{U}:=\left\{z\in\mathbb{C}^n:\operatorname{Im}z_n>|z^{\prime}|^2\right\}.
$$
where
$$
z=(z',z_n),\text{ where }z'=(z_1,\cdots,z_{n-1})\in\mathbb{C}^{n-1}\text{ and }z_n\in\mathbb{C}^1.
$$
Usually,for $p>0$,the space $L^{p}_{\alpha}(\mathcal{U},dV_\alpha)$ consists of all Lebesgue functions $f$ on $\mathcal{U}$ for which 
$$
\|f\|_{L^p_{\alpha}}:=\left\{\int_\mathcal{U}|f(z)|^pdV_\alpha(z)\right\}^{1/p}<\infty,
$$
Let $dV$ denote the Lebesgue measure on $\mathbb C^n$,for real parameter $\alpha>-1$,we define $dV_{\alpha}(z)=\boldsymbol\rho(z)^\alpha dV(z)$,where $\boldsymbol\rho(z)=\operatorname{Im}z_{n}-|z^{\prime}|^{2}. $

For given $a,b,c\in \mathbb R$,the following Forelli-Rudin type operators induced by liu et.al.\cite{ buttt},are defined by
$$
T_{a,b,c}f(z):=\boldsymbol{\rho}(z)^a\int_\mathcal{U}\frac{\boldsymbol{\rho}(w)^b}{\boldsymbol{\rho}(z,w)^c}f(w)dV(w),
$$
and
$$
S_{a,b,c}f(z):=\boldsymbol{\rho}(z)^a\int_{\mathcal{U}}\frac{\boldsymbol{\rho}(w)^b}{|\boldsymbol{\rho}(z,w)|^c}f(w)dV(w),
$$
where
$$
\boldsymbol{\rho}(z,w):=\frac{i}{2}(\bar{w}_n-z_n)-\langle z',w'\rangle .
$$

Recall the history of the boundedness over the Siegel upper half-space.In\cite{ButYur19} and \cite{Nor},the authors discussed the $L^p-L^q$ boundedness on the upper half-plane with specific parameter.Chen et.al.\cite{ButYur19} present the boundedness of $T_{0,0,c}$ from $L^p$ to $L^q$ for $(p,q)\in [1,\infty]\times[1,\infty]$.Sehba\cite{Nor} gave the boundedness of $T_{0,\alpha,2+\alpha-\gamma}$ from $L^p_\alpha$ to $L^q_\alpha$ for $1<p\leq q<\infty$,where $\gamma$ satisfies $0\leq \gamma <2+\alpha,\alpha>-1$.what's more,in \cite{but} Liu et.al extended the result of Chen et.al. to Siegel upper half-space and discussed the boundedness of $T_{0,0,c}$ from $L^p$ to $L^q$ space.The next year,in \cite{buttt},Liu et.al. gave the $L^p-L^p$ boundedness of $T_{a,b,c}$ and $S_{a,b,c}$,after that,two autuor of this paper also discussed the sufficiency and the necessity for the $L^p_\alpha-L^q_\beta$ boundedness of the two generalized classes of Forelli-Rudin type operators $T_{a,b,c}$ and $S_{a,b,c}$ in \cite{Pik91} for $1\leq p\leq q\leq \infty$.Finally,in 2023,Zhou et.al.\cite{butc} answered the Conjucture 1.4 of \cite{Pik91} for the case $p=1,q=\infty$ and established the relationship of the $L^p-L^q$ boundedness of Forelli-Rudin type operators between the unit ball $\mathbb B_n$ and the Siegel upper half-space $\mathcal{U}$.

Now,we extend the operators $T_{a,b,c}$ and $S_{a,b,c}$ to the multiparameter cases in this article. To be exact,we focus on the following two classes of multiparameter operators.
$$
T_{\vec{a},\vec{b},\vec{c}}f(z,w):=\boldsymbol{\rho}(z)^{a_1}\boldsymbol{\rho}(w)^{a_2}\int_{\mathcal{U}}\int_{\mathcal{U}}\frac{\boldsymbol{\rho}(u)^{b_1}\boldsymbol{\rho}(\eta)^{b_2}}{\boldsymbol{\rho}(z,u)^{c_1}\boldsymbol{\rho}(w,\eta)^{c_2}}f(u,\eta)dV(u)dV(\eta)
$$
and
$$
S_{\vec{a},\vec{b},\vec{c}}f(z,w):=\boldsymbol{\rho}(z)^{a_1}\boldsymbol{\rho}(w)^{a_2}\int_{\mathcal{U}}\int_{\mathcal{U}}\frac{\boldsymbol{\rho}(u)^{b_1}\boldsymbol{\rho}(\eta)^{b_2}}{|\boldsymbol{\rho}(z,u)|^{c_1}|\boldsymbol{\rho}(w,\eta)|^{c_2}}f(u,\eta)dV(u)dV(\eta)
$$
where $\vec{a}:=(a_{1},a_{2}),\vec{b}:=(b_{1},b_{2}),\vec{c}:=(c_{1},c_{2})\in\mathbb{R}^{2}.$

Obviously,these two integral operators are narural extension of the operators $T_{a,b,c}$ and $S_{a,b,c}$.In addition,our framwork for analyzing the boundedness of the Forlli-Rudin operators $T{\vec{a},\vec{b},\vec{c}}$ and $S_{\vec{a},\vec{b},\vec{c}}$ has also become weighted mixed-norm Lebesgue spaces.Notice that,for any given $\vec{p}:=(p_1,p_2)\in[1,\infty]^2$ and $\vec{\alpha}:=(\alpha_1,\alpha_2)\in(-1,\infty)^2$, the weighted mixed-norm Lebesgue space $L_{\vec{\alpha}}^{\vec{p}}:=L^{\vec{p}}(\mathcal{U}\times\mathcal{U},dV_{\alpha_1}dV_{\alpha_2})$ is defined to be the set of all measurable functions $f$ satisfy
$$
\|f\|_{L_{\vec{\alpha}}^{\vec{p}}}:=\left\{\int_{\mathcal{U}}\left[\int_{\mathcal{U}}|f(z,w)|^{p_1}dv_{\alpha_1}(z)\right]^{\frac{p_2}{p_1}}dv_{\alpha_2}(w)\right\}^{\frac{1}{p_2}}<\infty .
$$
with the usual modifications made when $p_i=\infty$,for some $i\in\{1,2\}$\{(see )\}.

With the notion of mixed-norm Lebesgue spaces, our main results are stated as follows. In the rest of this article,we always suppose $\vec{\alpha}:=(\alpha_{1},\alpha_{2})\in(-1,\infty)^{2}\mathrm{~and~}\vec{\beta}:=(\beta_{1},\beta_{2})\in(-1,\infty)^{2}.$

\begin{theorem}
Let $\vec{p}:=(p_1,p_2)$ and $\vec{q}:=(q_1,q_2)$ satisfy $1<p_-\leq p_+\leq q_-<\infty,\vec{\alpha}\neq \vec{\beta}$, where $p_{+}:=\max\{p_{1},p_{2}\},p_{-}:=\min\{p_{1},p_{2}\},$ and $  q_{-}:=\min\{q_{1},q_{2}\}.$ Then the following conditions are equivalent.\\
$(\mathrm{i})$The operator $T_{\vec{a},\vec{b},\vec{c}}$ is bounded from $L_{\vec{\alpha}}^{\vec{p}}$ to $L_{\vec{\beta}}^{\vec{q}}$.\\
$(\mathrm{ii})$The operator $S_{\vec{a},\vec{b},\vec{c}}$ is bounded from $L_{\vec{\alpha}}^{\vec{p}}$ to $L_{\vec{\beta}}^{\vec{q}}$.\\
$(\mathrm{iii})$The parameters satisfy that,for any $i\in\{1,2\}$,
\begin{equation}
\begin{cases}
 -q_ia_i<\beta_i+1,\alpha_i+1<p_i(b_i+1),\\
c_i\geq n+1+a_i+b_i+\frac{n+1+\beta_i}{q_i}-\frac{n+1+\alpha_i}{p_i}.   
\end{cases}   
\end{equation}
\end{theorem}

We also establish the following three $L_{\vec{\alpha}}^{\vec{p}}-L_{\vec{\beta}}^{\vec{q}}$ boundedness characterizations for the endpoints of exponents $\vec{p}:=(1,p_{2}),\vec{p}:=(p_{1},1)$,and $\vec{p}:=(1,1)$,respectively, where $1<p_{1},p_{2}<\infty.$ In addition,when $p_i=1$,we do not consider the case of $\alpha_i=b_i$.

\begin{theorem}
Let $\vec{p}:=(1,p_2)$ and $\vec{q}:=(q_1,q_2)$ satisfy $1<p_2\leq q_-<\infty,\vec{\alpha}\neq \vec{\beta}$, Then the following conditions are equivalent.\\
$(\mathrm{i})$The operator $T_{\vec{a},\vec{b},\vec{c}}$ is bounded from $L_{\vec{\alpha}}^{\vec{p}}$ to $L_{\vec{\beta}}^{\vec{q}}$.\\
$(\mathrm{ii})$The operator $S_{\vec{a},\vec{b},\vec{c}}$ is bounded from $L_{\vec{\alpha}}^{\vec{p}}$ to $L_{\vec{\beta}}^{\vec{q}}$.\\
$(\mathrm{iii})$The parameters satisfy that, for any $i\in\{1,2\}$,
\begin{equation}
 \begin{cases}
-q_ia_i<\beta_i+1,\alpha_1<b_1,\alpha_2+1<p_2(b_2+1),\\
c_1\geq a_1+b_1-\alpha_1+\frac{n+1+\beta_1}{q_1},\\
c_2\geq n+1+a_2+b_2+\frac{n+1+\beta_2}{q_2}-\frac{n+1+\alpha_2}{p_2}.
\end{cases}   
\end{equation}
\end{theorem}

\begin{theorem}
Let $\vec{p}:=(p_1,1)$ and $\vec{q}:=(q_1,q_2)$ satisfy $1<p_2\leq q_-<\infty,\vec{\alpha}\neq \vec{\beta}$, Then the following conditions are equivalent.\\
$(\mathrm{i})$The operator $T_{\vec{a},\vec{b},\vec{c}}$ is bounded from $L_{\vec{\alpha}}^{\vec{p}}$ to $L_{\vec{\beta}}^{\vec{q}}$.\\
$(\mathrm{ii})$The operator $S_{\vec{a},\vec{b},\vec{c}}$ is bounded from $L_{\vec{\alpha}}^{\vec{p}}$ to $L_{\vec{\beta}}^{\vec{q}}$.\\
$(\mathrm{iii})$The parameters satisfy that,for any $i\in\{1,2\}$,
\begin{equation}
\begin{cases}
-q_ia_i<\beta_i+1,\alpha_1+1<p_1(b_1+1),\alpha_2<b_2,\\
c_1\geq n+1+a_1+b_1+\frac{n+1+\beta_1}{q_1}-\frac{n+1+\alpha_1}{p_2},\\
c_2\geq a_2+b_2-\alpha_2+\frac{n+1+\beta_2}{q_2}.\\
\end{cases}    
\end{equation}
\end{theorem}

\begin{theorem}
Let $\vec{p}:=(1,1)$ and $\vec{q}:=(q_1,q_2)\in [1,\infty]$, $\vec{\alpha}\neq \vec{\beta}$, Then the following conditions are equivalent.\\
$(\mathrm{i})$The operator $T_{\vec{a},\vec{b},\vec{c}}$ is bounded from $L_{\vec{\alpha}}^{\vec{p}}$ to $L_{\vec{\beta}}^{\vec{q}}$.\\
$(\mathrm{ii})$The operator $S_{\vec{a},\vec{b},\vec{c}}$ is bounded from $L_{\vec{\alpha}}^{\vec{p}}$ to $L_{\vec{\beta}}^{\vec{q}}$.\\
$(\mathrm{iii})$The parameters satisfy that,for any $i\in\{1,2\}$,
\begin{equation}
\begin{cases}
-q_ia_i<\beta_i+1,\alpha_i<b_i,\\
c_i\geq a_i+b_i-\alpha_i+\frac{n+1+\beta_i}{q_i}.\\
\end{cases}    
\end{equation}
\end{theorem}

The rest of this article is organized as follows.

In Section 2,we show the proof of the necessary of boundedness for the operators $T_{\vec{a},\vec{b},\vec{c}}$.According to the case of the unit ball,the Key Lemma and Lemma 5 of \cite{buttt} are the most crusial tool for studying the boundedness of operators $T_{\vec{a},\vec{b},\vec{c}}$. we need to find suitable experimental functions $f_{\theta,\delta}$, and based on the definition of operator boundedness, we can find the conditions that the parameters $(\vec{a}, \vec{b}, \vec{c})$ in Theorem 1.1, 1.2, 1.3, and 1.4 satisfy. additionally, we also discussed the necessary conditions for operator $T_{\vec{a}, \vec{b}, \vec{c}}$ to be bounded when $(p_i,q_i)=(\infty,\infty),i\in\{1.2\}$.

In Section 3,we adopt the classic approach, which requires the use of Schur'test (see Proposition 3.1,Proposition 3.2,Proposition 3.3 and Proposition 3.4 of \cite{BondYur18}), which is generalized and applicable to mixed-norm Lebesgue spaces and proving the sufficiency of boundedness for $S_{0,\vec{b},\vec{c}}$,Adding observation 4.1 in \cite{BondYur18},we can obtain the sufficiency of boundedness for $S_{\vec{a},\vec{b},\vec{c}}$.In addition,on the details of further proof, we have solved the problem of “$\geq $”and found four specific parameters.After that,we conclude the work of section 2 and section 3 and gave the proof of Theorem 1.1-Theorem 1.4.

In section 4,we present two examples to illustrate the use of our main result.

Finally,due to the large number of constants generated during the proof process, in order to simplify the calculation, we also need to introduce some notation.Let C be a positive number,if $f\leq Cg$,then we write $f\lesssim g$,and if $f=Cg$,we write $f\sim g$.For any $p\in [1,\infty]$,We use $p'$ to represent its H$\ddot{o}$lder conjugate number, namely,$\frac{1}{p}+\frac{1}{p'}=1$.Similarly,For any $\vec{p}=(p_1,p_2)$,we write $\vec{p'}=(p'_1,p'_2)$,$\frac{1}{p_1}+\frac{1}{p'_1}=1$,$\frac{1}{p_2}+\frac{1}{p'_2}=1$.

\section{Proof of necessary for boundedness of $T_{\vec{a},\vec{b},\vec{c}}$}
In this section, we first introduce some important lemmas, and then study the necessary conditions for the boundedness of operators $T_ {\vec{a}, \vec{b}, \vec{c}}$.

\begin{lemma}[\cite{buttt},Key Lemma]
Let $s,t \in\mathbb{R}.$For any $z\in \mathcal{U}$,we have
\begin{equation}
\int\limits_{\mathcal{U}}\frac{\boldsymbol{\rho}(w)^t}{|\boldsymbol{\rho}(z,w)|^s}dV(w)=
\begin{cases}
\frac{C_2(n,s,t)}{\boldsymbol{\rho}(z)^{s-t-n-1}},&if\,t>-1\, and\, s-t>n+1\\+\infty,&otherwise
\end{cases}    
\end{equation}
where
$$
C_{2}(n,s,t):=\frac{4\pi^{n}\Gamma(1+t)\Gamma(s-t-n-1)}{\Gamma^{2}\left(s/2\right)}. 
$$
\end{lemma}

\begin{lemma}[\cite{buttt},Lemma 5]
assume $r,s>0,t>-1\,and\,r+s-t>n+1$,for any $z,u\in\mathcal{U}$,we have
\begin{equation}
\int_\mathcal{U}\frac{\boldsymbol{\rho}(w)^t}{\boldsymbol{\rho}(z,w)^r\boldsymbol{\rho}(w,u)^s}dV(w)=\frac{C_1(n,r,s,t)}{\boldsymbol{\rho}(z,u)^{r+s-t-n-1}}    
\end{equation}
where
$$
C_1(n,r,s,t):=\frac{4\pi^n\Gamma(1+t)\Gamma(r+s-t-n-1)}{\Gamma(r)\Gamma(s)}.
$$
\end{lemma}

\begin{lemma}
Let $\vec{p}:=(p_1,p_2)\in[1,\infty]^2,\vec{q}:=(q_1,q_2)\in[1,\infty]^2$,if $T_{\vec{a},\vec{b},\vec{c}}$ is bounded from $L_{\Vec{\alpha}}^{\Vec{p}}$ to $L_{\Vec{\beta}}^{\Vec{q}}$,then its adjoint operator $T^*_{\vec{a},\vec{b},\vec{c}}$ is bounded from $L_{\Vec{\beta}}^{\Vec{q'}}$ to $L_{\Vec{\alpha}}^{\Vec{p'}}$. $T^*_{\vec{a},\vec{b},\vec{c}}$ defined by setting
\[
\begin{split}
T^*_{\vec{a},\vec{b},\vec{c}}g(z,w):=\boldsymbol{\rho}(z)^{b_1-\alpha_1}\boldsymbol{\rho}(w)^{b_2-\alpha_2}\int_\mathcal{U}\int_\mathcal{U}\frac{\boldsymbol{\rho}(u)^{a_1+\beta_1}\boldsymbol{\rho}(\eta)^{a_2+\beta_2}}{\boldsymbol{\rho}(z,u)^{c_1}\boldsymbol{\rho}(w,\eta)^{c_2}}g(u,\eta)dV(u)dV(\eta).    
\end{split}
\]
\end{lemma}
    
\textbf{Proof:} Let $g\in L_{\Vec{\beta}}^{\Vec{q'}}$.From [\cite{DB22},304, Theorem 2]we have :
\[
\begin{split}
\|T_{\vec{a},\vec{b},\vec{c}}^*g\|_{L_{\vec{\alpha}}^{\vec{p}^{\prime}}}& =\sup_{\|f\|_{L_{\vec{\alpha}}^{\vec{p}}}=1}\left|\int_{\mathcal{U}}\int_{\mathcal{U}}\overline{f(z,w)}T_{\vec{a},\vec{b},\vec{c}}^{*}g(z,w)dV_{\alpha_{1}}(z)dV_{\alpha_{2}}(w)\right|  \\
&=\sup_{\left\|f\right\|_{L_{\vec{\alpha}}^{\vec{p}}}=1}\bigg|\int_{\mathcal{U}}\int_{\mathcal{U}}\overline{f(z,w)}\boldsymbol{\rho}(z)^{b_1-\alpha_1}\boldsymbol{\rho}(w)^{b_2-\alpha_2} \\
&\times\int_{\mathcal{U}}\int_{\mathcal{U}}\frac{\boldsymbol{\rho}(u)^{a_1+\beta_1}\boldsymbol{\rho}(\eta)^{a_2+\beta_2}}{\boldsymbol{\rho}(z,u)^{c_1}\boldsymbol{\rho}(w,\eta)^{c_2}}g(u,\eta)dV(u)dV(\eta)dV_{\alpha_{1}}(z)dV_{\alpha_{2}}(w)\bigg|.\\
\end{split}
\]
Further based on Fubini Theorem,\text{Hölder} Inequality
and $T_{\vec{a},\vec{b},\vec{c}}$ is bounded from $L_{\vec{\alpha}}^{\vec{p}}$ to $L_{\vec{\beta}}^{\vec{q}}$,we have 
\[
\begin{split}
\|T_{\vec{a},\vec{b},\vec{c}}^{*}g\|_{L_{\vec{\alpha}}^{\vec{p}^{\prime}}} =&\operatorname*{sup}_{\left\|f\right\|_{L_{\vec{\alpha}}^{\vec{p}}}=1}\left|\int_{\mathcal{U}}\int_{\mathcal{U}}g(u,\eta)\overline{{T_{\vec{a},\vec{b},\vec{c}}f(u,\eta)}}dV_{\beta_{1}}(u)dV_{\beta_{2}}(\eta)\right| \\
&\leq\sup_{\left\|f\right\|_{L_{\vec{\alpha}}^{\vec{p}}}=1}\left\|T_{\vec{a},\vec{b},\vec{c}}f\right\|_{L_{\vec{\beta}}^{\vec{q}}}\left\|g\right\|_{L_{\vec{\beta}}^{\vec{q}^{\prime}}}\lesssim\left\|g\right\|_{L_{\vec{\beta}}^{\vec{q}^{\prime}}}. 
\end{split}
\]

This finishes the proof.$\hfill\square$

\begin{lemma}
Let $1\leq p_{1}\leq q_{1}<\infty,1\leq p_{2}\leq q_{2}<\infty.\vec{\alpha}\neq \Vec{\beta}$,if $T_{\vec{a},\vec{b},\vec{c}}$ is boundedness from  $L^{\vec{p}}_{\Vec{\alpha}}$ to $L^{\vec{q}}_{\Vec{\beta}}$,then For any $i\in\{1,2\}$ we have 
$$
\begin{cases}-q_ia_i<\beta_i+1,\alpha_i+1<p_i(b_i+1),\\c_i\geq n+1+a_i+b_i+\frac{n+1+\beta_i}{q_i}-\frac{n+1+\alpha_i}{p_i}.\end{cases}
$$  
\end{lemma}

\textbf{Proof:} For any $\theta,\delta>0$,let
$$f_{\theta,\delta}(z,w):=\frac{\boldsymbol\rho(z)^t}{\rho(z,\theta\mathbf{i})^s}\frac{\boldsymbol\rho(w)^t}{\rho(z,\delta\mathbf{i})^s},\quad z,w\in\mathcal{U},$$
where $s,t$ is a real parameter and meets the following conditions
$$
\begin{aligned}
&s>0,\\
&t>\max\left\{-\frac{1+\alpha_i}{p_i},\:-1-b_i\right\},\\
&s-t>\max\left\{\frac{n+1+\alpha_i}{p_i}
,\:\frac{n+1+\beta_i}{q_i}+n+1+b_i-c_i\right\}.
\end{aligned}
$$
From Lemma 2.1 and Fubini Theorem we obtain
\[
\begin{split}
\|f_{\theta,\delta}\|_{L^{\vec{p}}_{\vec{\alpha}}}&=\left\{\int_{\mathcal{U}}\left[\int_{\mathcal{U}}|f_{\theta,\delta}(z,w)|^{p_1}dv_{\alpha_1}(z)\right]^{\frac{p_2}{p_1}}dv_{\alpha_2}(w)\right\}^{\frac{1}{p_2}}\\
&=\left[\int_{\mathcal{U}}\frac{\boldsymbol{\rho}(z)^{p_1t+\alpha_1}}{|\boldsymbol{\rho}(z,\theta\mathbf{i})|^{sp_1}}dV(z)\right]^{\frac{1}{p_1}}\left[\int_{\mathcal{U}}\frac{\boldsymbol{\rho}(w)^{p_2t+\alpha_2}}{|\boldsymbol{\rho}(w,\delta\mathbf{i})|^{sp_2}}dV(w)\right]^{\frac{1}{p_2}}\\
&\sim \theta^{\frac{n+1+\alpha_1}{p_1}+t-s}\cdot
\delta^{\frac{n+1+\alpha_2}{p_2}+t-s}.
\end{split}
\]
therefore $f_{\theta,\delta}\in L^{\vec{p}}(\mathcal{U}\times\mathcal{U},dV_{\alpha_{1}}dV_{\alpha_{2}}).$From Lemma 2.2,then 
\[
\begin{split}
&(T_{\vec{a},\vec{b},\vec{c}}f_{\theta,\delta})(z,w)\\
&=\boldsymbol{\rho}(z)^{a_1}\boldsymbol{\rho}(w)^{a_2}\int_{\mathcal{U}}\int_{\mathcal{U}}\frac{\boldsymbol{\rho}(u)^{b_1}\boldsymbol{\rho}(u)^t}{\boldsymbol{\rho}(z,u)^{c_1}\boldsymbol{\rho}(u,\theta\bold{i})^{s}}\frac{\boldsymbol{\rho}(\eta)^{b_2}\boldsymbol{\rho}(\eta)^t}{\boldsymbol{\rho}(w,\eta)^{s}\boldsymbol{\rho}(\eta,\delta\bold{i})^{s}}dV(u)dV(\eta)\\
&\sim\frac{\boldsymbol{\rho}(z)^{a_1}}{\boldsymbol{\rho}(z,\theta\bold{i})^{n+1+t+b_1-s-c_1}}\cdot 
\frac{\boldsymbol{\rho}(w)^{a_2}}{\boldsymbol{\rho}(w,\delta\bold{i})^{n+1+t+b_2-s-c_2}}
\end{split}
\]
Due to $T_{\vec{a},\vec{b},\vec{c}}f_{\theta,\delta}\in\ L^{\vec{q}}_{\vec{\beta}},$ according to Lemma 2.3,we have 
\begin{equation}
q_i a_i+\beta_i>-1,    
\end{equation}
\begin{equation}
q_i(c_i-a_i-b_i-n-1)+q_i(s-t)-n-1-\beta_i>0.   
\end{equation}
In addition,we write $c_i+s-b_1-t-n-1:=A_i$ also have
\[
\begin{split}
&\|T_{\vec{a},\vec{b},\vec{c}}f_{\theta,\delta}\|_{L^{\vec{q}}_{\beta}}\\
&\sim \left[\int_{\mathcal{U}}\frac{\boldsymbol{\rho}(z)^{q_1 a_1+\beta_1}}{|\boldsymbol{\rho}(z,\theta\mathbf{i})|^{q_1A_1}}dV(z)\right]^{\frac{1}{q_1}}\left[\int_{\mathcal{U}}\frac{\boldsymbol{\rho}(w)^{q_2 a_2+\beta_2}}{|\boldsymbol{\rho}(w,\beta\mathbf{i})|^{q_1A_2}}dV(w)\right]^{\frac{1}{q_2}}\\
&\sim \theta^{\frac{n+1+\beta_1}{q_1}-A_1} \cdot
\delta^{\frac{n+1+\beta_2}{q_2}-A_2}
\end{split}
\]
Due to $T_{\vec{a},\vec{b},\vec{c}}$ is bounded from $L^{\vec{p}}_{\Vec{\alpha}}$ to $L^{\vec{q}}_{\Vec{\beta}}$,then
\begin{equation}
\|T_{\vec{a},\vec{b},\vec{c}}f_{\theta,\delta}\|_{L_{\beta}^{\vec{q}}}\lesssim \|f_{\theta,\delta}\|_{L_{\alpha}^{\vec{p}}}.    
\end{equation}
Notice $\theta,\delta>0$,To make (2.5) valid, then
\begin{equation}
\frac{n+1+\beta_i}{q_i}-c_i-s+a_i+b_i+t+n+1\leq \frac{n+1+\alpha_i}{p_i}+t-s    
\end{equation}
For any $i\in\{1,2\}$,organize (2.6) to obtain 
$$
c_i\geq n+1+a_i+b_i+\frac{n+1+\beta_i}{q_i}-\frac{n+1+\alpha_i}{p_i}.
$$
Finally,From Lemma 2.3 and (2.3) we have 
$$
-p'_{i}(b_i-\alpha_i)<\alpha_i+1,
$$
When $p_i=1,p'_{i}=\infty$,$\alpha_i\leq b_i$,This proof only considers the case of $"<"$,then we also have
$$
\alpha_i+1<p_i(b_i+1).
$$
we complete the proof.$\hfill\square$

\textbf{Note:} When $\theta=\delta$, it can be inferred from equation (2.6) that for parameter $c_i$, we can also obtain the following weaker conclusion.
$$
\sum\limits_{i=1}^{2}c_i\geq \sum\limits_{i=1}^{2}n+1+a_i+b_i+\frac{n+1+\beta_i}{q_i}-\frac{n+1+\alpha_i}{p_i}.
$$

Let's consider the case where $p_i=1,\alpha_i=b_i $, and we can also obtain the following three lemmas.
\begin{lemma}
Let $1= p_{1}\leq q_{1}<\infty,1< p_{2}\leq q_{2}<\infty.\vec{\alpha}\neq \Vec{\beta}$,if $T_{\vec{a},\vec{b},\vec{c}}$ is bounded from $L^{\vec{p}}_{\Vec{\alpha}}$ to $L^{\vec{q}}_{\Vec{\beta}}$,then for any $i\in\{1,2\}$
$$
\begin{cases}-q_ia_i<\beta_i+1,\alpha_1=b_1,\alpha_2+1<p_2(b_2+1),\\
c_1= a_1+\frac{n+1+\beta_1}{q_1},\\
c_2\geq n+1+a_2+b_2+\frac{n+1+\beta_2}{q_2}-\frac{n+1+\alpha_2}{p_2}.
\end{cases}
$$
\end{lemma}

\begin{lemma}
Let $1<p_{1}\leq q_{1}<\infty,1= p_{2}\leq q_{2}<\infty.\vec{\alpha}\neq \Vec{\beta}$,if $T_{\vec{a},\vec{b},\vec{c}}$ is bounded from $L^{\vec{p}}_{\Vec{\alpha}}$ to $L^{\vec{q}}_{\Vec{\beta}}$,then for any $i\in\{1,2\}$
$$
\begin{cases}
-q_ia_i<\beta_i+1,\alpha_1+1<p_1(b_1+1),\alpha_2=b_2,\\
c_1\geq n+1+a_1+b_1+\frac{n+1+\beta_1}{q_1}-\frac{n+1+\alpha_1}{p_1},\\
c_2= a_2+\frac{n+1+\beta_2}{q_2}.
\end{cases}
$$
\end{lemma}

\begin{lemma}
Let $1= p_{1}\leq q_{1}<\infty,1= p_{2}\leq q_{2}<\infty.\vec{\alpha}\neq \Vec{\beta}$,if $T_{\vec{a},\vec{b},\vec{c}}$ is bounded from $L^{\vec{p}}_{\Vec{\alpha}}$ to $L^{\vec{q}}_{\Vec{\beta}}$,then for any $i\in\{1,2\}$
$$
\begin{cases}
-q_ia_i<\beta_i+1,\alpha_i=b_i,\\
c_i= a_i+\frac{n+1+\beta_i}{q_i}.\end{cases}
$$
\end{lemma}

We further analyze equation (2.6) and find that when $\frac{n+1+\beta_i}{q_i}=\frac{n+1+\alpha_i}{p_i} $, the term on the right side of the inequality is completely eliminated by the term on the left. At this point, $c_i=n+1+a_i+b_i $. Therefore, we can also obtain the following three lemmas (When $p_i=1$,without considering the case of $\alpha_i=b_i $).

\begin{lemma}
Let $1\leq p_{1}= q_{1}<\infty,1\leq p_{2}< q_{2}<\infty.\alpha_1=\beta_1$,if $T_{\vec{a},\vec{b},\vec{c}}$ is bounded from $L^{\vec{p}}_{\Vec{\alpha}}$ to $L^{\vec{q}}_{\Vec{\beta}}$,then for any $i\in\{1,2\}$
$$
\begin{cases}
-q_ia_i<\beta_i+1,\alpha_i+1<p_i(b_i+1),\\
c_1=n+1+a_1+b_1,\\
c_2\geq n+1+a_2+b_2+\frac{n+1+\beta_2}{q_2}-\frac{n+1+\alpha_2}{p_2}.\end{cases}
$$
\end{lemma}

\begin{lemma}
Let $1\leq p_{1}< q_{1}<\infty,1\leq p_{2}= q_{2}<\infty.\alpha_2=\beta_2$,if $T_{\vec{a},\vec{b},\vec{c}}$ is bounded from $L^{\vec{p}}_{\Vec{\alpha}}$ to $L^{\vec{q}}_{\Vec{\beta}}$,then for any $i\in\{1,2\}$
$$
\begin{cases}
-q_ia_i<\beta_i+1,\alpha_i+1<p_i(b_i+1),\\
c_1\geq n+1+a_1+b_1+\frac{n+1+\beta_1}{q_1}-\frac{n+1+\alpha_1}{p_1},\\
c_2=n+1+a_2+b_2.\\
\end{cases}
$$
\end{lemma}

\begin{lemma}
Let $\Vec{p}=\vec{q}\in [1,+\infty)^2,\vec{{\alpha}}=\vec{\beta}$,if $T_{\vec{a},\vec{b},\vec{c}}$ is bounded from $L^{\vec{p}}_{\Vec{\alpha}}$ to $L^{\vec{q}}_{\Vec{\beta}}$,then for any $i\in\{1,2\}$
$$
\begin{cases}
-q_ia_i<\beta_i+1,\alpha_i+1<p_i(b_i+1),\\
c_i=n+1+a_i+b_i.
\end{cases}
$$
\end{lemma}

Continuing, let's consider the case where $p_i=\infty $. Because according to previous research on classical Lebesgue spaces, we must have $q_i=\infty $.

\begin{lemma}
Let $p_{1}, q_{1}=\infty,1\leq p_{2}\leq q_{2}<\infty.\alpha_2\neq \beta_2$,if $T_{\vec{a},\vec{b},\vec{c}}$ is bounded from $L^{\vec{p}}_{\Vec{\alpha}}$ to $L^{\vec{q}}_{\Vec{\beta}}$,then
$$
\begin{cases}
a_1>0,b_1>-1,c_1=n+1+a_1+b_1,\\
-q_2a_2<\beta_2+1,\alpha_2+1<p_2(b_2+1),\\
c_2\geq n+1+a_2+b_2+\frac{n+1+\beta_2}{q_2}-\frac{n+1+\alpha_2}{p_2},
\end{cases}
$$\\
\end{lemma}

\textbf{Proof:} For any $\delta>0$,let
$$f_{\delta}(z,w):=\frac{\boldsymbol\rho(z,u)^{c_1}}{|\boldsymbol\rho(z,u)|^{c_1}}\frac{\boldsymbol\rho(w)^t}{\boldsymbol\rho(z,\delta\mathbf{i})^s},\quad z,w\in\mathcal{U},$$
where $s,t$ is a real parameter and meets the following conditions
$$
\begin{aligned}
&s>0,\\
&t>\max\left\{-\frac{1+\alpha_i}{p_i},\:-1-b_i\right\},\\
&s-t>\max\left\{\frac{n+1+\alpha_i}{p_i}
,\:\frac{n+1+\beta_i}{q_i}+n+1+b_i-c_i\right\}.
\end{aligned}
$$
From Lemma 2.1 and Fubini Theorem,we have
\[
\begin{split}
\|f_{\delta}\|_{L^{(\infty,p_2)}_{\alpha_2}}&=\left\{\int_{\mathcal{U}}\left[ess\sup_{z\in \mathcal{U}}|f_{\delta}(z,w)|\right]^{p_2}dv_{\alpha_2}(w)\right\}^{\frac{1}{p_2}}\\
&=\left[\int_{\mathcal{U}}\frac{\boldsymbol{\rho}(w)^{p_2t+\alpha_2}}{|\boldsymbol{\rho}(w,\delta\mathbf{i})|^{sp_2}}dV(w)\right]^{\frac{1}{p_2}}\\
&\sim\delta^{\frac{n+1+\alpha_2}{p_2}+t-s}
\end{split}
\]
Therefore $f_{\delta}\in L^{(\infty,p_2)}_{\alpha_2}(\mathcal{U}\times\mathcal{U}).$According to Lemma 2.2 we can obtain
\[
\begin{split}
&(T_{\vec{a},\vec{b},\vec{c}}f_{\delta})(z,w)\\
&\sim \boldsymbol{\rho}(z)^{a_1}\boldsymbol{\rho}(w)^{a_2}\int_{\mathcal{U}}\int_{\mathcal{U}}\frac{\boldsymbol{\rho}(u)^{b_1}}{|\boldsymbol{\rho}(z,u)|^{c_1}}\frac{\boldsymbol{\rho}(\eta)^{b_2}\boldsymbol{\rho}(\eta)^t}{\boldsymbol{\rho}(w,\eta)^{s}\boldsymbol{\rho}(\eta,\delta\bold{i})^{s}}dV(u)dV(\eta)\\
&\sim \boldsymbol{\rho}(z)^{(n+1+a_1+b_1-c_1)}\cdot 
\frac{\boldsymbol{\rho}(w)^{a_2}}{\boldsymbol{\rho}(w,\delta\bold{i})^{n+1+t+b_2-s-c_2}}
\end{split}
\]
Due to $T_{\vec{a},\vec{b},\vec{c}}f_{\delta}\in\ L^{\infty,q_2}_{\beta_2},$From Lemma2.2 we can see
$$q_2 a_2+\beta_2>-1,b_1>-1,a_1>0$$
$$q_2(c_2-a_2-b_2-n-1)+q_2(s-t)-n-1-\beta_2>0.$$
$$c_1=n+1+a_1+b_1$$
Further information can be obtained
\[
\begin{split}
\|T_{\vec{a},\vec{b},\vec{c}}f_{\delta}\|_{L^{(\infty ,q_2)}_{\beta_2}}&\sim \left[\int_{\mathcal{U}}\frac{\boldsymbol{\rho}(w)^{q_2 a_2+\beta_2}}{|\boldsymbol{\rho}(w,\delta\mathbf{i})|^{q_2A_2}}dV(w)\right]^{\frac{1}{q_2}}\\
&\sim \delta^{\frac{n+1+\beta_2}{q_2}-A_2}
\end{split}
\]
Due to $T_{\vec{a},\vec{b},\vec{c}}$ is bounded from $L^{(\infty,p_2)}_{\alpha_2}$ to $L^{(\infty,q_2)}_{\beta_2}$,we have
\begin{equation}
\|T_{\vec{a},\vec{b},\vec{c}}f_{\delta}\|_{L_{\beta}^{\vec{q}}}\lesssim\|f_{\delta}\|_{L_{\alpha}^{\vec{p}}}.    
\end{equation}
Because $\delta>0$,To make (2.7) valid,we have
\begin{equation}
\frac{n+1+\beta_2}{q_2}-c_2-s+a_2+b_2+t+n+1\leq \frac{n+1+\alpha_2}{p_2}+t-s     
\end{equation}
Organize (2.7) to obtain $c_2\geq n+1+a_2+b_2+\frac{n+1+\beta_2}{q_2}-\frac{n+1+\alpha_2}{p_2}.$
Finally,From Lemma2.3 and $-q_2a_2<\beta_2+1$ we can obtain
$$
-p'_{2}(b_2-\alpha_i)<\alpha_2+1,
$$
When $p_2=1,p'_{2}=\infty$,$\alpha_2\leq b_2$,This proof also only considers the case of $"<"$,then we also have
$$
\alpha_2+1<p_2(b_2+1).
$$
we complete the proof.$\hfill\square$\\

According to the proof of Lemma 2.11, we can also obtain the following two lemmas.\\

\begin{lemma}
Let $p_{2}, q_{2}=\infty,1\leq p_{1}\leq q_{1}<\infty.\alpha_1\neq \beta_1$,if $T_{\vec{a},\vec{b},\vec{c}}$ is bounded from $L^{\vec{p}}_{\Vec{\alpha}}$ to $L^{\vec{q}}_{\Vec{\beta}}$,then
$$
\begin{cases}
a_2>0,b_2>-1,c_2=n+1+a_2+b_2,\\
-q_1a_1<\beta_1+1,\alpha_1+1<p_1(b_1+1),\\
c_1\geq n+1+a_1+b_1+\frac{n+1+\beta_1}{q_1}-\frac{n+1+\alpha_1}{p_1},
\end{cases}
$$\\
\end{lemma}

\textbf{Proof:} For any $\theta>0$,let
$$f_{\theta}(z,w):=\frac{\boldsymbol\rho(z)^{t}}{\boldsymbol\rho(z,\theta\mathbf{i})^{s}}\frac{\boldsymbol\rho(w,\eta)^{c_2}}{|\boldsymbol\rho(w,\eta)|^{c_2}},\quad z,w\in\mathcal{U},$$
Referring to the proof of Lemma 2.11, Lemma 2.12 can be obtained.$\hfill\square$\\

\begin{lemma}
Let $\Vec{p},\vec{q}=(\infty,\infty)$,if $T_{\vec{a},\vec{b},\vec{c}}$ is bounded from $L^{\vec{p}}_{\Vec{\alpha}}$ to $L^{\vec{q}}_{\Vec{\beta}}$,then
$$
\begin{cases}
a_i>0,b_i>-1,\\
c_i=n+1+a_i+b_i,\\
\end{cases}\qquad i\in \{1,2\}
$$\\
\end{lemma}

\textbf{Proof:} We let 
$$f(z,w):=\frac{\boldsymbol\rho(z,u)^{c_1}}{|\boldsymbol\rho(z,u)|^{c_1}}\frac{\boldsymbol\rho(w,\eta)^{c_2}}{|\boldsymbol\rho(w,\eta)|^{c_2}},\quad z,w\in\mathcal{U},$$
Referring to the proof of Lemma 2.11, Lemma 2.13 can be obtained.$\hfill\square$\\

\section{Proof of sufficiency for boundedness of $S_{{\vec{a},\vec{b},\vec{c}}}$}

This section aims to show the sufficiency of the boundedness of multiparameter Forelli-Rudin type operators $S_{\vec{a},\vec{b},\Vec{c}}$.The key is the use of Schur'test(see Proposition 3.1,Proposition 3.2,Proposition 3.3 and Proposition 3.4 of \cite{BondYur18}), which is promoted for mixed-norm Lebesgue spaces,With the help of these Schur’s tests, we next prove the sufficiency of the main theorems.\\

\begin{lemma}
Let $1=p_-<p_+\leq q_-<\infty,1= p_+ \leq q_- < \infty$ or $1< p_- \leq p_+ \leq q_- < \infty , \vec{\alpha}\neq \vec{\beta}.$ For any 
$i\in\{1,2\}$,
$$
\begin{cases}-q_ia_i<\beta_i+1,\alpha_i+1<p_i(b_i+1),\\c_i\geq n+1+a_i+b_i+\frac{n+1+\beta_i}{q_i}-\frac{n+1+\alpha_i}{p_i},\end{cases}
$$\\
then $S_{\vec{a},\vec{b},\vec{c}}$ is bounded from $L^{\vec{p}}_{\Vec{\alpha}}$ to $L^{\vec{q}}_{\Vec{\beta}}$.
\end{lemma}

\textbf{Approach:} Let's first consider the boundedness of $S_{0, \vec {b}, \vec {c} }$, and then obtain the boundedness of $S_{\vec {a}, \vec {b}, \vec {c}} $ based on Observation 4.2 in \cite{BondYur18}.\\

\textbf{Proof:} For any $i\in \{1,2\}$,\,let
$$\lambda_i:=\frac{n+1+\beta_i}{q_i}-\frac{n+1+\alpha_i}{p_i}$$
$$c_i:=n+1+b_i+\lambda_i+\omega_i,\tau_i:=c_i-b_i-\omega_i+\alpha_i.$$
where $\omega_i\geq 0.$ Due to $-\frac{1+\beta_i}{q_i}<0$,we can know that there exist two negative numbers $r_1,r_2,$ such that,for any $i\in\{1,2\},$$-\frac{1+\beta_i}{q_i}<r_i<0.$ Additionally,from $\alpha_i+1<p_i(b_i+1),$ we can obtain $\alpha_i+1<p_i(b_i+\omega_i+1),$ write $b_i+\omega_i=b'_i$,it,follows that,for any $i\in\{1,2\},$
\begin{equation}
 b'_i-\alpha_i+\frac{\alpha_i+1}{p'_i}>0.   
\end{equation}
and 
\begin{equation}
\tau_i=\frac{n+1+\alpha_i}{p_i^{\prime}}+\frac{n+1+\beta_i}{q_i}>0,  
\end{equation}
combined with (3.1),further presents that
$$
-\frac{\tau_i(1+\alpha_i)}{p_i'}-\frac{(b'_i-\alpha_i)(n+1+\alpha_i)}{p_i'}<\frac{(b'_i-\alpha_i)(n+1+\beta_i)}{q_i}.
$$
Then there exist $s_1,s_2$ such that,for any $i\in\{1,2\},$
$$
-\frac{\tau_i(1+\alpha_i)}{p_i'}-\frac{(b'_i-\alpha_i)(n+1+\alpha_i)}{p_i'}<\tau_is_i+(b'_i-\alpha_i)(s_i-r_i)<\frac{(b'_i-\alpha_i)(n+1+\beta_i)}{q_i},
$$
Organize and obtain
\begin{equation}
-\frac{1+\alpha_i}{p_i'}-(b'_i-\alpha_i)\gamma_i<s_i<(b'_i-\alpha_i)\delta_i,   
\end{equation}
where, for any $i\in\{1,2\},$
$$
\gamma_{i}:=\frac{(n+1+\alpha_{i})/p_{i}^{\prime}+s_{i}-r_{i}}{\tau_{i}},\delta_{i}:=\frac{(n+1+\beta_{i})/q_{i}+r_{i}-s_{i}}{\tau_{i}}.
$$
Obviously,$\gamma_i+\delta_i=1.$Below we provide the test function that satisfies Schur'test
$$
h_{1}(u,\eta):=\boldsymbol{\rho}(u)^{s_{1}}\boldsymbol{\rho}(\eta
)^{s_{2}},h_{2}(z,w):=\boldsymbol{\rho}(z)^{r_{1}}\boldsymbol{\rho}(w)^{r_{2}},
$$
$$
K_{1}(z,u):=\frac{\boldsymbol{\rho}(u)^{b'_1-\alpha_1}}{|\boldsymbol{\rho}(z,u)|^{c_1}},K_{2}(w,\eta):=\frac{\boldsymbol{\rho}(\eta
)^{b'_2-\alpha_2}}{|\boldsymbol{\rho}(w,\eta)|^{c_2}}.
$$
Therefore,
$$
S_{0,\vec{b},\vec{c}}f(z,w)=\int_{\mathcal{U}}\int_{\mathcal{U}}K_{1}(z,u)K_{2}(w,\eta)f(u,\eta)dV_{\alpha_{1}}(u)dV_{\alpha_{2}}(\eta).
$$
We next prove this lemma by considering the following two cases.\\

$\boldsymbol{ Case 1.}$ $p_->1$,this article also found specific $r_i,s_i,i\in\{1,2\}$ as follows
$$
r_i=\frac{-(1+\beta_i)}{p_1p'_1q_1q'_1},\qquad s_i=\frac{n}{p_1p'_1q_1q'_1}.
$$
We use Proposition 3.1 in \cite{BondYur18} to obtain the boundedness of $S_ {0, \vec {b}, \Vec {c} }$, we need to consider
\[
\begin{split}
 \int_{\mathcal{U}}&[K_1(z,u)]^{\gamma_1p_1'}[K_2(w,\eta)]^{\gamma_2p_1'}[h_1(u,\eta)]^{p_1'}dV_{\alpha_1}(u) \\
 &=\frac{\boldsymbol\rho(\eta)^{(b'_2-\alpha_2)p'_1\gamma_2+s_2p'_1}}{|\boldsymbol{\rho}(w,\eta)|^{c_2\gamma_2p'_1}}\int_{\mathcal{U}}\frac{\boldsymbol\rho(u)^{(b'_1-\alpha_1)p'_1\gamma_1+s_1p'_1+\alpha_1}}{|\boldsymbol{\rho}(z,u)|^{c_1\gamma_1p'_1}}dV(u).
\end{split}
\]
From the left inequality of (3.3),we know
\begin{equation}
(b'_i-\alpha_i)\gamma_ip_i'+s_ip_i'+\alpha_i>-1.    
\end{equation}
In addition,by the fact that $(c_i-b'_1+\alpha_i)=\tau_i\gamma_i,$it follows that,for any $i\in\{1,2\}$,
\begin{equation}
c_i\gamma_ip'_i-n-1-(b'_i-\alpha_i)\gamma_ip'_i-\alpha_i-s_ip'_i=-r_ip'_i>0.    
\end{equation}
From (3.4),(3.5) and Lemma 2.1,we have
$$
\int_{\mathcal{U}}\frac{\boldsymbol\rho(u)^{(b'_1-\alpha_1)p'_1\gamma_1+s_1p'_1+\alpha_1}}{|\boldsymbol{\rho}(z,u)|^{c_1\gamma_1p'_1}}dV(u)\sim \boldsymbol\rho(z)^{r_1p'_1},
$$
Similarly,we can obtain
\[
\begin{split}
\int_{\mathcal{U}}\bigg[\int_{\mathcal{U}}&[K_1(z,u)]^{\gamma_1p_1'}[K_2(w,\eta)]^{\gamma_2p_1'}[h_1(u,\eta)]^{p_1'}dV_{\alpha_1}(u)\bigg]^{p_2'/p_1'}dV_{\alpha_2}(\eta)\\
&\sim\boldsymbol\rho(z)^{r_1p'_2}\int_{\mathcal{U}}\frac{\boldsymbol\rho(\eta)^{(b'_2-\alpha_2)p'_2\gamma_2+s_2p'_2+\alpha_2}}{|\boldsymbol{\rho}(w,\eta)|^{c_2\gamma_2p'_2}}dV(\eta)\\
&\sim\boldsymbol\rho(z)^{r_1p'_2}\boldsymbol\rho(w)^{r_2p'_2}\sim[h_2(z,w)]^{p'_2}.
\end{split}
\]
Thus,Satisfying (3.1) in Proposition 3.1 in\cite{BondYur18}, we next will continue to check (3.2) in Proposition 3.1 in\cite{BondYur18}.Observe that
\[
\begin{split}
 \int_{\mathcal{U}}&[K_1(z,u)]^{\delta_1q_1}[K_2(w,\eta)]^{\delta_2q_1}[h_2(z,w)]^{q_1}dV_{\beta_1}(z) \\
 &=\frac{\boldsymbol\rho(\eta)^{(b'_2-\alpha_2)\delta_2q_1
 }\boldsymbol\rho(w)^{r_2q_1}}{|\boldsymbol{\rho}(w,\eta)|^{c_2\delta_2q_1}}\int_{\mathcal{U}}\frac{\boldsymbol\rho(u)^{(b'_1-\alpha_1)\delta_1q_1}\boldsymbol\rho(z)^{r_1q_1+\beta_1}}{|\boldsymbol{\rho}(z,u)|^{c_1\delta_1q_1}}dV(z).
\end{split}
\]
According to the definition of $r_i $, it can be inferred that
\begin{equation}
r_iq_i+\beta_i>-1.    
\end{equation}
Moreover,for any $i\in\{1,2\}$,$(c_i-b'_i+\alpha_i)\delta_i=\tau_i\delta_i=\frac{n+1+\beta_i}{q_i}+r_i-s_i,$ combined with (3.3) will obtain
\begin{equation}
 c_i\delta_iq_i-n-1-r_iq_i-\beta_i=(b'_i-\alpha_i)\delta_iq_i-s_iq_i>0.   
\end{equation}
From (3.6),(3.7) and Lemma 2.1,we have
$$
\int_{\mathcal{U}}\frac{\boldsymbol\rho(z)^{r_1q_1+\beta_1}}{|\boldsymbol{\rho}(z,u)|^{c_1\delta_1q_1}}dV(z)\sim \boldsymbol\rho(u)^{s_1q_1-(b'_1-\alpha_1)\delta_1q_1}.
$$
Similarly,we have
\[
\begin{split}
\int_{\mathcal{U}}&\left[\int_{\mathcal{U}}[K_1(z,u)]^{\delta_1q_1}[K_2(w,\eta)]^{\delta_2q_1}[h_2(z,w)]^{q_1}dV_{\beta_1}(z)\right]^{q_2/q_1}dV_{\beta_2}(w)\\
&\sim\boldsymbol\rho(u)^{s_1q_1}\boldsymbol\rho(u)^{(b'_2-\alpha_2)\delta_2q_2}\int_{\mathcal{U}}\frac{\boldsymbol\rho(w)^{r_2q_2+\beta_2}}{|\boldsymbol{\rho}(w,\eta)|^{c_2\delta_2q_2}}dV(w)\\
&\sim \boldsymbol\rho(u)^{s_1q_1}\boldsymbol\rho(\eta)^{s_2q_2}.
\end{split}
\]
Hence,we complete the proof of Lemma 3.1 in this case.$\hfill\square$\\

$\boldsymbol{ Case 2.}$ $p_+=1$.For any $i\in\{1,2\},$According to the proof process above, note that $\gamma_{i}:=\frac{s_{i}-r_{i}}{\tau_{i}},\delta_{i}:=\frac{(n+1+\beta_{i})/q_{i}+r_{i}-s_{i}}{\tau_{i}},$ and 
\begin{equation}
-(b'_i-\alpha_1)\gamma_i<s_i<(b'_i-\alpha_i)\delta_i.    
\end{equation}

We use Proposition 3.2 in \cite{BondYur18} to obtain the boundedness of $S_ {0, \vec {b}, \Vec {c} }$,  need to consider
$$
[K_1(z,u)]^{\gamma_1}[K_2(w,\eta)]^{\gamma_2}h_1(u,\eta)=\frac{\boldsymbol{\rho}(u)^{\gamma_1(b'_1-\alpha_1)+s_1}\boldsymbol{\rho}(\eta)^{\gamma_2(b'_2-\alpha_2)+s_2}}{|\boldsymbol{\rho}(z,u)|^{c_1\gamma_1}|\boldsymbol{\rho}(w,\eta)|^{c_2\gamma_2}}.
$$
From (\cite{DV19},Lemma 2.1),it follows that,for any $z,w\in\mathcal{U}$
$$
2|\boldsymbol{\rho}(z,w)|\geq\max\{\boldsymbol{\rho}(z),\boldsymbol{\rho}(w)\}.
$$
and
\begin{equation}
\frac{\boldsymbol{\rho}(z)}{|\boldsymbol{\rho}(z,u)|}\leq \frac{1}{2},\qquad \frac{\boldsymbol{\rho}(u)}{|\boldsymbol{\rho}(z,u)|}\leq \frac{1}{2}.   
\end{equation}
Since $\tau_i=c_i-b'_i+\alpha_i$\, and\, $\gamma_i=(s_i-r_i)/\tau_i$,we have
\begin{equation}
c_i\gamma_i=(b'_i-\alpha_i)\gamma_i+s_i-r_i.    
\end{equation}
By the left inequality of (3.8),we know that,for any $i\in\{1,2\},\gamma_i(b'_i-\alpha_i)+s_i>0,$\,and the fact that $-r_i>0$.Combined with (3.9) (3.10),for any given $z\in\mathcal{U}$,for any $u\in\mathcal{U},$we have
\begin{equation}
\frac{\boldsymbol{\rho}(u)^{\gamma_1(b'_1-\alpha_1)+s_1}\boldsymbol{\rho}(z)^{-r_1}}{|\boldsymbol{\rho}(z,u)|^{c_1\gamma_1}}=\left(\frac{\boldsymbol{\rho}(u)}{|\boldsymbol{\rho}(z,u)|}\right)^{\gamma_1(b'_1-\alpha_1)+s_1}\left(\frac{\boldsymbol{\rho}(z)}{|\boldsymbol{\rho}(z,u)|}\right)^{-r_1}\lesssim 1.   
\end{equation}
Similarly,for any given $w\in\mathcal{U}$,for any $\eta\in\mathcal{U},$we have
$$
\frac{\boldsymbol{\rho}(\eta)^{\gamma_2(b'_2-\alpha_2)+s_2}\boldsymbol{\rho}(w)^{-r_2}}{|\boldsymbol{\rho}(w,\eta)|^{c_2\gamma_2}}\lesssim 1.
$$
This, together with (3.20), further conclude that, for any given $(z,w)\in \mathcal{U}\times \mathcal{U}$,
$$
\mathop{}_{(u,\eta)\in\mathcal{U\times \mathcal{U}}}^{\text{ess sup}}[K_{1}(z,u)]^{\gamma_{1}}[K_{2}(w,\eta)]^{\gamma_{2}}h_{1}(u,\eta)\lesssim h_{2}(z,w).
$$
Therefore, the test for (3.3) in Proposition 3.2 in \cite{BondYur18} has been completed.Notice that  (3.4) is the same as (3.2) in Proposition 3.1, because there is no change in $q_i $. Hence,we complete the proof of Lemma 3.1 in this case.$\hfill\square$\\

$\boldsymbol{ Case 3.}$ $1=p_1<p_2\leq q_- $. According to the proof process above,note that $\gamma_{1}:=\frac{s_{i}-r_{i}}{\tau_{i}},\gamma_{2}:=\frac{(n+1+\alpha_2)/p'_2+s_{i}-r_{i}}{\tau_{i}},\delta_{i}:=\frac{(n+1+\beta_{i})/q_{i}+r_{i}-s_{i}}{\tau_{i}}$,and

\begin{equation}
-(b'_1-\alpha_1)\gamma_1<s_1<(b'_1-\alpha_1)\delta_1.    
\end{equation}
\begin{equation}
 -\frac{1+\alpha_2}{p_2'}-(b'_2-\alpha_2)\gamma_2<s_2<(b'_2-\alpha_2)\delta_2   
\end{equation}

We use Proposition 3.3 in \cite{BondYur18} to obtain the boundedness of $S_ {0, \vec {b}, \Vec {c} }$,we need to check (3.5) in \cite{BondYur18}.with the proof of $\boldsymbol{ Case 2}$,we know that,for any $u\in\mathcal{U}$
\begin{equation}
\begin{split}
[K_1(z,u)]^{\gamma_1}[K_2(w,\eta)]^{\gamma_2}h_1(u,\eta)&=\frac{\boldsymbol{\rho}(u)^{\gamma_1(b'_1-\alpha_1)+s_1}\boldsymbol{\rho}(\eta)^{\gamma_2(b'_2-\alpha_2)+s_2}}{|\boldsymbol{\rho}(z,u)|^{c_1\gamma_1}|\boldsymbol{\rho}(w,\eta)|^{c_2\gamma_2}}\\   &\lesssim \boldsymbol{\rho}(z)^{-r_1}\frac{\boldsymbol{\rho}(\eta)^{\gamma_2(b'_2-\alpha_2)+s_2}}{|\boldsymbol{\rho}(w,\eta)|^{c_2\gamma_2}}
\end{split}    
\end{equation}

Further from the proof of $\boldsymbol{ Case 1}$ and (3.14),obtain
$$
\int_{\mathcal{U}}\left[\mathop{}_{u\in\mathcal{U}}^{\text{ess sup}}[K_{1}(z,u)]^{\gamma_{1}}[K_{2}(w,\eta)]^{\gamma_{2}}h_{1}(u,\eta)\right]^{p'_2}dV_{\alpha_2}(\eta)\lesssim\left[h_{2}(z,w)\right]^{p'_2}.
$$
Therefore, the test for (3.5) in Proposition 3.3 in \cite{BondYur18} has been completed.Notice that  (3.6) is the same as (3.2) in Proposition 3.1, because there is no change in $q_i $. Hence,we complete the proof of Lemma 3.1 in this case.$\hfill\square$\\

$\boldsymbol{ Case 4.}$ $1=p_2<p_1\leq q_-.$This situation involves swapping $p_1 $\,and\, $p_2 $in $\boldsymbol {Case 3} $, and then verifying Proposition 3.4 in \cite{BondYur18}. The process is extremely similar, so this process is omitted.Hence,we have completed the proof of Lemma 3.1.$\hfill\square$\\

According to the proof of $\boldsymbol{ Case 4.}$, we can also obtain the following three lemmas.

\begin{lemma}
Let $1< p_1=q_1< \infty ,1<p_2<q_1\leq q_2<\infty ,\alpha_1= \beta_1,$ for any 
$i\in\{1,2\}$,
$$
\begin{cases}-q_ia_i<\beta_i+1,\alpha_i+1<p_i(b_i+1),\\
c_1=n+1+a_1+b_1,\\
c_2\geq n+1+a_2+b_2+\frac{n+1+\beta_2}{q_2}-\frac{n+1+\alpha_2}{p_2},\\
\end{cases}
$$\\
then $S_{\vec{a},\vec{b},\vec{c}}$ is bounded from $L^{\vec{p}}_{\Vec{\alpha}}$ to $L^{\vec{q}}_{\Vec{\beta}}$.
\end{lemma}

\begin{lemma}
Let $1< p_2=q_2< \infty ,1<p_1<q_2\leq q_1<\infty ,\alpha_2= \beta_2,$ for any 
$i\in\{1,2\}$,
$$
\begin{cases}-q_ia_i<\beta_i+1,\alpha_i+1<p_i(b_i+1),\\
c_1\geq n+1+a_1+b_1+\frac{n+1+\beta_1}{q_1}-\frac{n+1+\alpha_1}{p_1},\\
c_2=n+1+a_2+b_2,\\
\end{cases}
$$\\
then $S_{\vec{a},\vec{b},\vec{c}}$ is bounded from $L^{\vec{p}}_{\Vec{\alpha}}$ to $L^{\vec{q}}_{\Vec{\beta}}$.
\end{lemma}

\begin{lemma}
Let $1<p_1=p_2=q_1=q_2<\infty$, for any 
$i\in\{1,2\}$,
$$
\begin{cases}-q_ia_i<\beta_i+1,\alpha_i+1<p_i(b_i+1),\\
c_i=n+1+a_i+b_i,\\
\end{cases}
$$\\
then $S_{\vec{a},\vec{b},\vec{c}}$ is bounded from $L^{\vec{p}}_{\Vec{\alpha}}$ to $L^{\vec{q}}_{\Vec{\beta}}$.
\end{lemma}

\textbf{Conclusion:} According to the Lemma 2.4 and Lemma 3.1,we can obtain the proof of Theorem 1.1-1.4.

\section{Applications}

Based on the main results we obtained above, we provide two corresponding theorems.

After that,we need to introduce more notation.We know that Bergman kernel function $K_{\Omega}$ induces a Riemannian metric on a domain $\Omega\in\mathbb{C}^n$.The infinitesimal Bergman metric is defined by
$$
g_{i,j}^{\Omega}(z)=\frac{1}{n+1}\frac{\partial^{2}\log K_{\Omega}(z,z)}{\partial z_{i}\partial\bar{z}_{j}},\quad i,j=1,2,\ldots,n,
$$
and the complex matrix
$$
B(z)=\left(g_{i,j}^{\Omega}(z)\right)_{1\leq i,j\leq n}.
$$
is called the Bergman matrix on $\Omega$.For a $\mathbb{C}^1$ curve $\gamma :[0,1]\to\Omega$,the Bergman
length of $\gamma$ is defined by
$$
\ell(\gamma):=\int_{0}^{1}\langle B(\gamma(t))\gamma'(t),\gamma'(t)\rangle dt.
$$
If $z,w\in\mathcal{U},$The Bergman distance between them is
$$
\delta_\Omega(z,w):=\inf\{\ell(\gamma):\gamma(0)=z,\gamma(1)=w\},
$$
where infimum is taken over all $\mathbb{C}^1$ curves from $z$ to $w$.If $\Omega_1,\Omega_2\in \mathbb{C}^n$,and $\psi$ is a biholomorphic mapping of $\Omega_1$ onto $\Omega_2$,then for any $z,w\in\Omega_1,\delta_{\Omega_1}(z,w)=\delta_{\Omega_2}(\psi(z),\psi(w)),$from (19) in \cite{buttt},we have
\begin{equation}
 \begin{split}
\delta_{\mathcal{U}}(z,w)&=\delta_{\mathbb{B}}(\Phi^{-1}(z),\Phi^{-1}(w))=\tanh^{-1}\left(\left|\varphi_{\Phi^{-1}(z)}(\Phi^{-1}(w)\right|\right)\\
&=\tanh^{-1}\sqrt{1-\frac{\rho(z)\rho(w)}{|\rho(z,w)|^{2}}}.
\end{split}   
\end{equation}
Next, we consider the following operators
\[
\begin{split}
&S_{\Vec{a},\vec{b},\Vec{c}}^{\vec{d}}f(z,w)\\
&=\boldsymbol{\rho}(z)^{a_1}\boldsymbol{\rho}(w)^{a_2}\int_{\mathcal{U}}\int_{\mathcal{U}}\frac{\boldsymbol{\rho}(u)^{b_1}\delta_{\mathcal{U}}(z,u)^{d_1}}{|\boldsymbol{\rho}(z,u)|^{c_1}}\frac{\boldsymbol{\rho}(\eta)^{b_2}\delta_{\mathcal{U}}(w,\eta)^{d_2}}{|\boldsymbol{\rho}(w,\eta)|^{c_2}}f(u,\eta)dV(u)dV(\eta).   
\end{split}
\]
where $\Vec{a},\Vec{b},\Vec{c},\vec{d}\in \mathbb{R}^2.$We can obtain the following theorem.\\

\begin{theorem}
Assume $\vec{\alpha},\vec{\beta}\in (1,\infty)^2,\vec{p},\vec{q}\in [1,\infty)^2,$ for any 
$i\in\{1,2\}$,
$$
\begin{cases}
-q_ia_i<\beta_i+1,\alpha_i+1<p_i(b_i+1),d_i\geq 0,\\
c_i> n+1+a_i+b_i+\frac{n+1+\beta_i}{q_i}-\frac{n+1+\alpha_i}{p_i},
\end{cases}
$$
then $S_{\vec{a},\vec{b},\vec{c}}^{\vec{d}}$ is bounded from $L^{\vec{p}}_{\Vec{\alpha}}$ to $L^{\vec{q}}_{\Vec{\beta}}$.    
\end{theorem}

\textbf{Proof:} Pick $\epsilon>0$ so small such that,
$$
\begin{cases}
-q_i(a_i-d_i\epsilon)<\beta_i+1,\alpha_i+1<p_i(b_i+1-d_i\epsilon),d_i\geq 0,\\
c_i-2d_i\epsilon> n+1+a_i+b_i+\frac{n+1+\beta_i}{q_i}-\frac{n+1+\alpha_i}{p_i}.
\end{cases}
$$
Since for any $x>0,\epsilon>0,\log x<x^{\epsilon},$from (4.1) we know
$$
\delta_{\mathcal{U}}(z,w)\lesssim\log\frac{4|\boldsymbol{\rho}(z,w)|^2}{\boldsymbol{\rho}(z)\boldsymbol{\rho}(w)}\lesssim1+\frac{|\boldsymbol{\rho}(z,w)|^{2\epsilon}}{\boldsymbol{\rho}(z)^\epsilon\boldsymbol{\rho}(w)^\epsilon}.
$$
we write
$$
\begin{cases}
\vec{a^{1}}=(a_1-d_1\epsilon,a_2)\\
\vec{b^{1}}=(b_1-d_1\epsilon,b_2)\\
\vec{c^{1}}=(c_1-2d_1\epsilon,c_2)
\end{cases}
\begin{cases}
\vec{a^{2}}=(a_1,a_2-d_2\epsilon)\\
\vec{b^{2}}=(b_1,b_2-d_2\epsilon)\\
\vec{c^{2}}=(c_1,c_2-2d_2\epsilon)
\end{cases}
\begin{cases}
\vec{a^{3}}=(a_1-d_1\epsilon,a_2-d_2\epsilon)\\
\vec{b^{3}}=(b_1-d_1\epsilon,b_2-d_2\epsilon)\\
\vec{c^{3}}=(c_1-2d_1\epsilon,c_2-2d_2\epsilon)
\end{cases}
$$
After calculation, it can be concluded that
\[
\begin{split}
|S_{\Vec{a},\vec{b},\Vec{c}}^{\vec{d}}f(z,w)|&\lesssim |S_{\Vec{a},\vec{b},\Vec{c}}(|f|)(z,w)|+|S_{\Vec{a^{1}},\vec{b^{1}},\Vec{c^{1}}}(|f|)(z,w)|\\
&+|S_{\Vec{a^{2}},\vec{b^{2}},\Vec{c^{2}}}(|f|)(z,w)|+|S_{\Vec{a^{3}},\vec{b^{3}},\Vec{c^{3}}}(|f|)(z,w)|
\end{split}
\]
From Lemma 3.1,we can know $S_{\vec{a},\vec{b},\vec{c}}^{\vec{d}}$ is bounded from $L^{\vec{p}}_{\Vec{\alpha}}$ to $L^{\vec{q}}_{\Vec{\beta}}$.$\hfill\square$\\

To prove Theorem 4.3, we need to use the following lemma.\\
\begin{lemma}
 If $\vec{\gamma}:=(\gamma_{1},\gamma_{2})\in(-1,\infty)^{2},\Vec{p}\in [1,\infty)^2,f\in A_{\Vec{\gamma}}^{\vec{p}},$ then for any $z,w\in \mathcal{U},$
$$
f(z,w)=\int_{\mathcal{U}}\int_{\mathcal{U}}\frac{f(u,\eta)}{\boldsymbol{\rho}(z,u)^{n+1+\gamma_1}\boldsymbol{\rho}(w,\eta)^{n+1+\gamma_2}}dV_{\gamma_1}(u)dV_{\gamma_2}(\eta).
$$    
\end{lemma}

\textbf{Proof} From Fubini Theorem and Theorem 2.1 in \cite{dju},we have
\[
\begin{split}
&\int_{\mathcal{U}}\int_{\mathcal{U}}\frac{f(u,\eta)}{\boldsymbol{\rho}(z,u)^{n+1+\gamma_1}\boldsymbol{\rho}(w,\eta)^{n+1+\gamma_2}}dV_{\gamma_1}(u)dV_{\gamma_2}(\eta)\\
&=\int_{\mathcal{U}}\frac{1}{\boldsymbol{\rho}(w,\eta)^{n+1+\gamma_2}}\int_{\mathcal{U}}\frac{f(u,\eta)}{\boldsymbol{\rho}(z,u)^{n+1+\gamma_1}}dV_{\gamma_1}(u)dV_{\gamma_2}(\eta)\\
&=\int_{\mathcal{U}}\frac{f(z,\eta)}{\boldsymbol{\rho}(w,\eta)^{n+1+\gamma_2}}dV_{\gamma_2}(\eta)\\
&=f(z,w)
\end{split}
\]
as asserted.$\hfill\square$\\

We denote by $A_{\vec{\alpha}}^{\vec
{p}}(\mathcal{U}\times \mathcal{U})$ the Bergman space,which is the closed subspace of $L^{\vec{p}}(\mathcal{U}\times\mathcal{U},\boldsymbol{\rho}^{{\alpha_1}}\boldsymbol{\rho}^{{\alpha_2}})$ consisting of holomorphic functions on $\mathcal{U}\times \mathcal{U}$,we write  $\partial_{n_z}:=\partial/(\partial z_n).$ The
following result plays an important role in the study of the Besov spaces over the Siegel upper half-space.\\

\begin{theorem}
Suppose $\vec{p},\vec{q}\in (1,\infty)^2,$ for any  $i\in\{1,2\},\alpha_i>\frac{1}{p_i}-1,\beta_i>\frac{1}{q_i}-1,N\in \mathbb{N}$,then $\partial_{n_z}^{N}\partial_{n_w}^{N}$ is a bounded
linear operator from $A_{\vec{\alpha}}^{\vec{p}}$ to $A_{{\vec{\beta}+N\vec{q}}}^{\vec{q}}$.\\    
\end{theorem}

\textbf{Proof:} From Lemma 4.2 and the suppose,we have
$$
f(z,w)=\int_{\mathcal{U}}\int_{\mathcal{U}}\frac{f(u,\eta)}{\boldsymbol{\rho}(z,u)^{n+1+\alpha_1}\boldsymbol{\rho}(w,\eta)^{n+1+\alpha_2}}dV_{\alpha_1}(u)dV_{\alpha_2}(\eta).
$$
we write
$$
\begin{cases}
\Vec{N}=(N,N)\\
\Vec{N_\beta}=(n+1+\beta_1+N,n+1+\beta_2+N)
\end{cases}
$$
It follows that
\[
\begin{split}
 &\left|\boldsymbol{\rho}(z)^N\boldsymbol{\rho}(w)^N\partial_{n_w}^N\partial_{n_z}^Nf(z,w)\right|\\
 &\lesssim\boldsymbol{\rho}(z)^N\boldsymbol{\rho}(w)^N\int_\mathcal{U}\int_\mathcal{U}|f(u,\eta)|\frac{\boldsymbol{\rho}(u)^{\beta_1}}{|\boldsymbol{\rho}(z,u)|^{n+1+\beta_1+N}}\frac{\boldsymbol{\rho}(\eta)^{\beta_2}}{|\boldsymbol{\rho}(w,\eta)|^{n+1+\beta_2+N}}dV(u)dV(\eta) \\
 &=S_{\vec{N},\vec{\beta},\vec{N_\beta}}(|f|)(z,w).
\end{split}
\]
According to Lemma 3.1,we obtain
$$
\left\|\partial_{n_{z}}^N\partial_{n_{w}}^Nf\right\|_{\Vec{\beta}+N\Vec{q}}^{\vec{q}}=\left\|\boldsymbol{\rho}^{N_z}\boldsymbol{\rho}^{N_w}\partial_{n_{z}}^N\partial_{n_{w}}^Nf\right\|_{\vec{\alpha}}^{\vec{p}}\lesssim\|f\|_{\alpha}^{\vec{p}}.
$$ 
as asserted.$\hfill\square$\\

\medskip
\noindent
{\bf Hongheng Yin}\\
School of Mathematical Sciences.\\
Beijing Normal University,China. \\
E-mail: YongHeng@mail.bnu.edu.cn

\medskip
\noindent
{\bf Guan-Tie Deng}\\
School of Mathematical Sciences.\\
Beijing Normal University,China. \\
E-mail: denggt@mail.bnu.edu.cn

\medskip
\noindent
{\bf Zhi-Qiang Gao}\\
School of Mathematical Sciences.\\
Beijing Normal University,China. \\
E-mail: gaozq@mail.bnu.edu.cn

\end{document}